\documentclass[A4]{article}

\DeclareSymbolFont{lettersA}{U}{txmia}{m}{it}
\DeclareMathSymbol{\FFF}{\mathord}{lettersA}{'206}
\DeclareMathSymbol{\NNN}{\mathord}{lettersA}{'216}
\DeclareMathSymbol{\RRR}{\mathord}{lettersA}{'222}
\DeclareMathSymbol{\ZZZ}{\mathord}{lettersA}{'232}
\DeclareMathSymbol{\QQQ}{\mathord}{lettersA}{'221}
\DeclareMathSymbol{\CCC}{\mathord}{lettersA}{'203}

\usepackage{color}

\usepackage[dvipdfmx]{graphicx}

\usepackage{amsmath,amsthm,amssymb}

\newcommand{\wsq}{\hfill$\square$}

\begin{document}

\title{EXISTENCE OF PERIODIC ORBITS FOR\\SINGULAR-HYPERBOLIC LYAPUNOV STABLE SETS}
\author{Kouta Nakai\\Graduate School of Mathematical Sciences, The University of Tokyo}
\date{2009/01/19}
\maketitle

\noindent ABSTRACT.
In \cite{BM}, Bautista and Morales proved the existence of periodic orbits in singular-hyperbolic attracting sets.
In this paper, we extend their result to singular-hyperbolic Lyapunov stable sets.

\section{INTRODUCTION}
	In 1963, E. N. Lorenz published the so-called Lorenz equation:
	\[
	\begin{array}{ll}
	x'=-10x+10y\\
	y'=28x-y-xz\\
	z'=xy-8/3z
	\end{array}
	\]
	that is related to some of the unpredictable behavior of the weather\cite{L}.
	Later in 1979 Guckenheimer and Williams\cite{GW} and
	in 1982 Afraimovich, Bykov and Shilnikov\cite{ABS}
	introduced the Geometric Lorenz Attractor(GLA).
        In 1999, Tucker\cite{T} showed that the GLA indeed corresponds to
        the behaviour of solutions of the original Lorenz equation.
        The GLA allowed us to examine qualitative behavior.
        It has been shown that the GLA has a sort of hyperbolicity and dense periodic orbits\cite{PM}.
	After that, singular-hyperbolicity was introduced as an extended concept of hyperbolicity.
        The GLA is an example of singular-hyperbolic attractor.
        By the existence of a transitive orbit and the Shadowing Lemma\cite{HK},
        it is known that every singular-hyperbolic attractor contains a periodic orbit.
        Then, it is natural to ask whether so does every singular hyperbolic attracting set or not.
        This problem was solved affirmatively by Bautista and Morales\cite{BM}.        
        However, it is not known whether every singular-hyperbolic Lyapunov stable set is attracting or not\cite{CM}.
	(It is known that every isolated Lyapunov stable set is attracting\cite{CM2}.)
        So, it is still worth proving that singular-hyperbolic Lyapunov stable sets has periodic orbits.\\

Let $M$ be a compact 3-manifold and let $X_t$ be a $C^1$ flow on $M$.
We denote the vector field associated to $X_t$ by $X$.
	Given $p\in M$, an $orbit$ of $X_t$ is the set $O_X(p)=\{ X_t(p);\; t\in \RRR \}$.
	In particular, \textit{positive orbit} means $\{ X_t(p); \; t\geq 0\}$.
	We denote the \textit{omega-limit set} and the \textit{alpha-limit set} of a point $p$
        by $\omega _X(p)$ and $\alpha _X(p)$ respectively.
A \textit{singularity} of $X_t$ is a point $\sigma \in M$ such that $X(\sigma )=0$.
	We denote the set of all singularities of $X_t$ by $\mathrm{Sing}(X)$, and
        singularities in a subset $B\subset M$ by $\mathrm{Sing}_X(B)$.
A \textit{periodic orbit} of $X_t$ is 
an orbit $O_X(p)$ such that $X_T(p)=p$ for some $p\in M$ and $T>0$.
A \textit{closed orbit} of $X_t$ is either a singularity or a periodic orbit of $X_t$.
A compact set $\Lambda \subset M$ is \textit{invariant} if $X_t(\Lambda )=\Lambda $ for all $t\in \RRR$.

\noindent
	\textbf{Definition 1.}
	A compact invariant set $\Lambda\subset M$ is \textit{Lyapunov stable} 
	if for given neighborhood $U$ of $\Lambda$, there is a neighborhood 
	$V$ of $\Lambda $ in $U$ such that the positive orbit of every point in $V$ is contained in $U$ i.e.,
	\[	X_t(x)\in U \; \mathrm{for}\; \mathrm{all}\; t\geq 0 \; \mathrm{and}\; \mathrm{for}\; \mathrm{all}\; x\in V.	\]
        We denote by $V_{Lyp}(U)$ such a neighborhood $V$.
\\

\noindent
\textbf{Definition 2.}
A compact invariant set $\Lambda$ of $X_t$ is \textit{hyperbolic} if there are positive constants $K$, $\lambda $ and a
continuous invariant tangent bundle decomposition $T_{\Lambda}M=E_{\Lambda}^s\oplus E_{\Lambda}^X\oplus E_{\Lambda}^u$ such that
\begin{enumerate}
	\item $E_{\Lambda}^s$ is \textit{contracting}, i.e.,
		\[ \Vert DX_t(x)\vert_{E_x^s} \Vert \leq Ke^{-\lambda t} \; \mathrm{for \; all} \; t>0 \; \mathrm{and\;for\;all}\; x\in \Lambda.\]
	\item $E_{\Lambda}^u$ is \textit{expanding}, i.e.,
		\[ \Vert DX_{-t}(x)\vert_{E_x^u} \Vert \leq Ke^{-\lambda t} \; \mathrm{for \; all} \; t>0 \; \mathrm{and\;for\;all}\; x\in \Lambda.\]
	\item $E_{\Lambda}^X$ is tangent to the vector field $X$ associated to $X_t$
\end{enumerate}

For a linear space or a submanifold $L$ of $M$ we denote the dimension of $L$ by $\mathrm{dim}(L)$.
	By the Invariant Manifold Theorem\cite{HPS}, for a hyperbolic set $\Lambda$ of $X_t$ and $p\in\Lambda$,
	the \textit{strong stable manifold} $W^{ss}_X(p)$ of $p$ and the \textit{strong unstable manifold} $W^{uu}_X(p)$ of $p$
	exist and they are $C^1$ submanifolds of $M$:
	\begin{eqnarray*}
	&& W_X^{ss}(p)=\{ x\in M;\; \lim_{t\rightarrow\infty}d(X_t(x),X_t(p))=0\} \\
	&& W_X^{uu}(p)=\{ x\in M;\; \lim_{t\rightarrow -\infty}d(X_t(x),X_t(p))=0\}
	\end{eqnarray*}

It is known that $W_X^{ss}(p)$ and $W_X^{uu}(p)$ are tangent respectively to the linear spaces $E_p^s$ and $E_p^u$ at $p$.
A closed orbit $O$ of $X_t$ is \textit{hyperbolic} if it is hyperbolic as a compact invariant set.
A hyperbolic closed orbit $O$ is \textit{saddle-type} if $E_p^s\neq 0$ and $E_p^u\neq 0$ for some (and hence for all) $p\in O$.\\
For a linear operator $A$, we denote the minimum norm by $m(A)=\inf_{v\neq 0}(\Vert Av \Vert / \Vert v \Vert )$\\

\noindent
\textbf{Definition 3.}
Let $\Lambda $ be a compact invariant set of $X_t$. A continuous invariant splitting $T_{\Lambda}M=E_{\Lambda}\oplus F_{\Lambda}$ 
over $\Lambda$ is \textit{dominated} if there are positive constants $K$ and $\lambda$ such that
 \[ \frac{\Vert DX_t(x)\vert_{E_x}\Vert}{m(DX_t(x)\vert_{F_x})}\leq Ke^{-\lambda t} \; \mathrm{for \; all}\; t>0 \; \mathrm{and \; for \; all}\; x\in \Lambda \]
Hereafter we assume that $E_x \neq 0$ and $F_x \neq 0$ for every $x \in \Lambda$. A compact invariant set $\Lambda$ 
is \textit{partially hyperbolic} if it exhibits a dominated splitting $T_{\Lambda}M=E_{\Lambda}^s\oplus E_{\Lambda}^c$ 
such that $E^s_{\Lambda}$ is \textit{contracting}, i.e.,
 \[ \Vert DX_t(x)\vert_{E_x^s} \Vert \leq Ke^{-\lambda t} \: \mathrm{for \: all} \: t>0 \: \mathrm{and\:for\:all}\: x\in \Lambda. \]

	Now we define singular-hyperbolicity.
\\

\noindent
\textbf{Definition 4.}
A \textit{singular-hyperbolic set} $\Lambda$ of $X_t$ is a partially hyperbolic set with a
\textit{volume expanding central subbundle} $E_{\Lambda}^c$, i.e.,
	\[ \vert\det(DX_t(x)\vert_{E_x^c})\vert \geq K^{-1}e^{\lambda t}\; \mathrm{for \; all}\; t>0 \; \mathrm{and \; for \; all}\; x\in \Lambda \]
and all of singularities in $\Lambda$ are hyperbolic.\\

	A \textit{singular-hyperbolic Lyapunov stable set} is a singular-hyperbolic set which is simultaneously Lyapunov stable.
	Similarly \textit{singular-hyperbolic attracting set} and \textit{singular-hyperbolic attractor} are defined.
	Here, a compact invariant set $\Lambda$ is an \textit{attracting} set if it has a positively invariant isolating block $U$
        (i.e., $\bigcap_{t\in\RRR}X_t(U)=\Lambda$ and $X_t(U)\subset U$ for $\forall t\geq 0$)
        and is an \textit{attractor} if it is a transitive
        (i.e., for $\forall U$, $V\subset \Lambda$ there exists $t\geq 0$ such that $X_t(U)\cap V \neq \emptyset$) attracting set.\\

Let $\Lambda$ be a non-trivial connected singular-hyperbolic set of $X_t$.
	Non-trivial means that it is not a closed orbit.
	For the singular-hyperbolic splitting $T_{\Lambda}M=E_{\Lambda}^s\oplus E_{\Lambda}^c$, it is known that
	$X(x)\subset E_x^c$ and $\mathrm{dim}(E_x^s)=1$ for any $x\in\Lambda$ [1,Theorem 3].
\\

\noindent
\textbf{Definition 5.}
A singularity is \textit{Lorenz-like} if it has real eigenvalues $\lambda_1$, $\lambda_2$ and $\lambda_3$ satisfing
	\[ 	\lambda_2<\lambda_3<0<-\lambda_3<\lambda_1. 	\]
	We denote the set of Lorenz-like singularities of $X_t$ in a subset $B\subset M$ by $\mathrm{LSing}_X(B)$.
\\

We introduce some invariant manifolds asociated to a Lorenz-like singularity $\sigma$.
Since $\sigma$ is hyperbolic, the stable and the unstable manifolds $W_X^s(\sigma)$ and $W_X^u(\sigma)$ exist.
They are tangent at $\sigma$ to the eigenspaces associated to the set of eigenvalues $\{ \lambda_2,\lambda_3\}$ and $\{ \lambda_1\}$ respectively.
In particular, $W_X^s(\sigma)$ is two-dimensional and $W_X^u(\sigma)$ is one-dimensional.
A further invariant manifold called the \textit{strongly stable manifold} $W^{ss}_X(\sigma)$
exists and is tangent at $\sigma$ to the eigenspace associated to the subset of eigenvalue $\{ \lambda_2\}$.
	For a singularity which has two positive eigenvalues, we put $W_X^{ss}(\sigma)=W_X^s(\sigma)$.
        The following property is known [1, Lemma 1]:
for a connected singular-hyperbolic set $\Lambda$ of $X_t$,
if $\sigma\in\mathrm{Sing}_X(\Lambda)$, then $\sigma$ is Lorenz-like or has two positive eigenvalues.
Moreover in any case we have that $\Lambda\cap W_X^{ss}(\sigma)=\{ \sigma \}$.\\

	Here we recall a few examples of singular-hyperbolic sets. 
        The GLA is an example of a singular-hyperbolic attractor with periodic orbits.
        An example of singular-hyperbolic attracting set with periocdic orbits was
        recently provided by Morales[11,Theorem B]
        (which is constructed by modifing the Cherry-flow\cite{PM} and the GLA.). 
        On the other hands, there exists an example of a singular-hyperbolic set without periodic orbits.
        It is a flow on a solid torus ($D^1\times S^1$) constructed by Morales\cite{M2} using the Cherry-flow.
\\

\noindent
\textbf{Theorem.}
Every singular-hyperbolic Lyapunov stable set of a $C^1$ flow has a periodic orbit.\\

        Let us give a brief sketch of the proof.
        Let $\Lambda$ be a singular-hyperbolic Lyapunov stable set.
        If there is a singularity in $\Lambda$, it is known that the singularity is a Lorenz-like or has two positive eigenvalues.
        We consider dividing into the following three cases.
        The case where there are no singularities in $\Lambda$,
        the case where there are singularities except for Lorenz-like ones,
        and the case where there is a Lorenz-like singularity.

        In the first case, $\Lambda$ is a (saddle-type) hyperbolic set.
        Take $x\in\Lambda$ and a cross-section $\Sigma$ with $x$.
        By the Shadowing Lemma\cite{HK}, there exists a periodic point $p$ near $x$ and hence near $\Lambda$ in $\Sigma$.
        Assume $p\notin\Lambda$ and choose $U\supset\Lambda$ with $p\notin U$.
        Then, take $V_{Lyp}(U)$.
        Since the stable and the unstable manifolds of $x$ and $p$ are large enough to intersect transversally,
        using the $\lambda$-lemma\cite{PM}, we can see that some image of any neighborhood of $x$
        contains a point arbitarily close to $p$, contradicting the Lyapunov stability.

        In the second case, we can show that, for any $x\in\Lambda\backslash\mathrm{Sing}_X(\Lambda)$,
        $\omega_X(x)$ does not contain singularities.
        Then, it is a (saddle-type) hyperbolic set.
        Take $y\in\omega_X(x)$ and a cross-section $\Sigma$ with $y$.
        Applying the Shadowing Lemma\cite{HK}, we have a periodic point $p$ near $y$ and hence near $\Lambda$ in $\Sigma$.
        Assume $p\notin\Lambda$ and choose $U\supset\Lambda$ with $p\notin U$.
        Then, take $V_{Lyp}(U)$.
        Since the stable and the unstable manifolds of $y$ and $p$ are large enough to intersect transversally,
        similarly to the first case, this contradicts the Lyapunov stability.

	In the last case, if there is a Lorenz-like singularity, we construct cross-sections near all of Lorenz-like ones.
        It is proved that the return map on these sections satisfies some conditions.
        When the return map satisfies these conditions, some iteration under the return map of any curve on the sections horizontally crosses
        an element of a finite set of vertical bands of the sections.
        Using this, we have some iteration of an element of the finite set horizontally crosses another element.
        Repeating this, we obtain a chain of vertical bands.
        By the finiteness of our vertical bands, we obtain a closed sub-chain, calling a cycle.
        There exists a periodic point determined by a cycle.
        Now, we take a sequence of curves $\{ c_n\}$ accumulating on $\Lambda$.
        By the argument above, we have a cycle $\beta_n$ coming from each curve $c_n$.
        Since cycles are also finite, there exists a cycle $\beta$ as an accumulation point of $\{ \beta_n\}$.
        Take a sub-sequence $\{ \tilde{c}_n\}$ corresponds to the cycle $\beta$.
        Let $p$ be a periodic point determined by $\beta$.
        Assume $p\notin\Lambda$, and take neighborhoods $U$ and $V_{Lyp}(U)$ of $\Lambda$ as above.
        Take $N$ with $\tilde{c}_N\subset V_{Lyp}(U)$, then the image of some iteration of $\tilde{c}_N$ horizontally crosses $W^s_X(p)$.
        This contradicts the Lyapunov stability again.
\\

	In the proof of the Theorem, we use many lemmas of \cite{BM},
        many of which are easy to extend to the Lyapunov stability condition.
        Thus we omit their proofs and give only statements.
        The extension of Lemma 2 has a difficulty, so we describe it in detail.
        For the proof of the Theorem,
        we have difficulties mainly in proving that a periodic point near $\Lambda$ is indeed contained in $\Lambda$.
\\

        In Section 2, we prepare some settings for the proof of the Theorem.
        Then, in Section 3 we prove the Theorem.
        In Appendix, we give an example of Lyapunov stable set which is not attracting
        with a close property to singular-hyperbolicity.

%%------------------------------------------------------------------------------------------------------------------------
%%------------------------------------------------------------------------------------------------------------------------

\section{PRELIMINARIES}
We consider certain maps called \textit{hyperbolic triangular maps}
defined on a finite disjoint union of copies of $[-1,1]\times [-1,1]$
and discontinuous maps still preserving a continuous vertical foliation.
We also assume two hypotheses (H1) and (H2)
imposing certain amount of differentiability close to the point whose iteration falls eventually in the interior of $\Lambda$.

Proposition 1 asserts the existence of a hyperbolic periodic point for the hyperbolic triangular map that satisfies
(H1) and (H2) and has the large domain. 
	Then, we construct a family of cross-sections, so-called the \textit{singular cross-section}.

\subsection{Hyperbolic triangular maps}

Let $I=[-1,1]$ be a unit closed interval. Let $I_i$ be a copy of $I$ and let $\Sigma_i$ be a copy of the square 
$I^2=I\times I$ for $i=1,2,\ldots ,k$. We denote the disjoint union of the squares $\Sigma_i$ by $\Sigma$.
Put
	\[ 	L_{-i}=\{ -1\} \times I_i, \quad L_{0i}=\{ 0\} \times I_i \quad \mathrm{and} \quad L_{+i}=\{ 1\} \times I_i	 \]
for $i=1,2,\ldots ,k$ and 
	\[ 	L_{-}=\bigcup_{i=1}^{k}L_{-i}, \quad L_{0}=\bigcup_{i=1}^{k}L_{0i} \quad \mathrm{and} \quad L_{+}=\bigcup_{i=1}^{k}L_{+i}.	 \]
Given a map $F$, we denote the domain of $F$ by Dom($F$).
	A point $x\in \mathrm{Dom}(F)$ is \textit{periodic} for $F$ if there is an integer $n\geq 1$ such that 
	$F^j(x)\in \mathrm{Dom}(F)$ for all $0\leq j\leq n-1$ and $F^n(x)=x$.
        We denote all the periodic points of $F$ by $\mathrm{Per}(F)$.

A \textit{curve} $c$ in $\Sigma$ is the image of a $C^1$ injective map $c:\mathrm{Dom}(c)\subset \RRR \rightarrow \Sigma$ 
with Dom($c$) being a compact interval. We often identify $c$ with its image set. A curve $c$ is \textit{vertical} 
if it is the graph of a $C^1$ map $g:I_i\rightarrow I_i$, i.e.,
	\[ 	c=\{ (g(y),y); \; y\in I_i \subset \Sigma_i \; \mathrm{for \; some}\; i=1,2,\ldots ,k\}. 	\]

	A continuous foliation $\mathcal{F}_i$ on a component $\Sigma_i$ is called \textit{vertical} if its leaves are vertical curves 
	and $L_{-i}$, $L_{0i}$ and $L_{+i}$ are also leaves of $\mathcal{F}_i$.
A \textit{vertical foliation} $\mathcal{F}$ of $\Sigma$ is a foliation which restricted to each component $\Sigma_i$ 
of $\Sigma$ is a vertical foliation.
It follows that the leaves $L$ of a vertical foliation $\mathcal{F}$ are vertical curves 
hence differentiable ones. In particular, the tangent space $T_xL$ is well defined for all $x\in L$.
	For a foliation $\mathcal{F}$, we use the notation $L\in\mathcal{F}$ to mean that $L$ is a leaf of $\mathcal{F}$.

	For a map $F:\mathrm{Dom}(F)\subset \Sigma \rightarrow \Sigma $ and a vertical foliation $\mathcal{F}$ on $\Sigma$,
we say that $F$ \textit{preserves} $\mathcal{F}$ if for every leaf $L$ of $\mathcal{F}$ containd in $\mathrm{Dom}(F)$, there is 
a leaf $f(L)$ of $\mathcal{F}$ suth that $F(L)\subset f(L)$ and the restriction to $L$, $F\vert_L:L\rightarrow f(L)$ is continuous.
        For a vertical foliation $\mathcal{F}$ on $\Sigma$, a subset $B\subset \Sigma$ is \textit{saturated set} for 
	$\mathcal{F}$ if $B$ is the union of leaves of $\mathcal{F}$. We say that $B$ is $\mathcal{F}$-saturated for short.
	For a subset A of $\Sigma$ denote by $\mathcal{F}_A$ the union of leaves that intersects A.
	If $A=\{ x\} $, then $\mathcal{F}_x$ is the leaf of $\mathcal{F}$ containing $x$.	
	For a subsets $A$, $B\subset \Sigma$, we say that \textit{A covers B} if $\mathcal{F}_A \supset B$.
\\

        Now we define the \textit{triangular map} and consider its hyperbolicity.
\\

\noindent
\textbf{Definition 6.}
A map $F:\mathrm{Dom}(F)\subset \Sigma \rightarrow \Sigma $ is called \textit{triangular} if it preserves a 
vertical foliation $\mathcal{F}$ on $\Sigma$ such that Dom($F$) is $\mathcal{F}$-saturated.\\

We define the hyperbolicity of triangular maps with cone fields in $\Sigma$.
We denote the tangent bundle of $\Sigma$ by $T\Sigma$.
Given $x\in\Sigma$, $\alpha>0$ and a linear subspace $V_x\subset T_x\Sigma$, we denote the cone around $V_x$ in $T_x\Sigma$ with inclination $\alpha$
by $C_{\alpha}(x,V_x)\equiv C_{\alpha}(x)$, namely
	\[ C_{\alpha}(x)=\{ v_x\in T_x\Sigma; \; \angle(v_x,V_x)\leq \alpha \}. \]
Here $\angle(v_x,V_x)$ denotes the angle between a vector $v_x$ and the subspace $V_x$.
A \textit{cone field} in $\Sigma$ is a continuous map $C_{\alpha}:x\in\sigma\rightarrow C_{\alpha}(x)\subset T_x\Sigma$,
where $C_{\alpha}(x)$ is a cone with constant inclination $\alpha$ on $T_x\Sigma$.
A cone field $C_{\alpha}$ is called \textit{transversal} to a vertical foliation $\mathcal{F}$ on $\Sigma$ if
$T_xL$ is not contained in $C_{\alpha}(x)$ for any $x\in L$ and $L\in \mathcal{F}$.\\

\noindent
\textbf{Definition 7.}
Let $F:\mathrm{Dom}(F)\subset \Sigma \rightarrow \Sigma$ be a triangular map with associated vertical foliation $\mathcal{F}$.
Given $\lambda>0$ we say that $F$ is $\lambda$-hyperbolic if there is a cone field $C_{\alpha}$ in $\Sigma$ such that
\begin{enumerate}
	\item 	$C_{\alpha}$ is transversal to $\mathcal{F}$.
	\item 	If $x\in \mathrm{Dom}(F)$ and $F$ is differentiable at $x$, then\\
		$DF(x)(C_{\alpha}(x))\subset \mathrm{Int}(C_{\frac{\alpha}{2}}(F(x)))$ and
                $\Vert DF(x)\cdot v_x\Vert \geq \lambda \cdot \Vert v_x\Vert$
                for all $v_x\in C_{\alpha}(x)$.
\end{enumerate}

\subsection{Hypotheses (H1) and (H2).}
They impose some regularity around those leaves whose iteration \textit{eventually fall into} $\Sigma\backslash (L_-\cup L_+)$.
To state them we need the following definitions.\\

\noindent
\textbf{Definition 8.}
Let $F:\mathrm{Dom}(F)\subset\Sigma\rightarrow\Sigma$ be a triangular map such that $L_-\cup L_+\subset \mathrm{Dom}(F)$.
For all $L\in\mathcal{F}$ contained in Dom($F$) we define the (possibly $\infty$) number $n(L)$ as follows:
\begin{enumerate}
	\item If $F(L)\subset\Sigma\backslash(L_-\cup L_+)$, we define $n(L)=0$.
	\item If $F(L)\subset L_-\cup L_+$, we define
        	\[ n(L)=\sup\{ n\geq 1;\;F^i(L)\subset \mathrm{Dom}(F) \; \mathrm{and} \; F^{i+1}(L)\subset L_-\cup L_+ \; \mathrm{whenever} \; 0\leq i\leq n-1\} .\]
\end{enumerate}

Essentially $n(L)+1$ gives the first non-negative iterate of $L$ falling into $\Sigma\backslash(L_-\cup L_+)$.\\

\noindent
\textbf{Definition 9.}
Let $F:\mathrm{Dom}(F)\subset \Sigma\rightarrow\Sigma$ be a triangular map such that $L_-\cup L_+\subset\mathrm{Dom}(F)$.
We say that $F$ satisfies:
\begin{description}
	\item[(H1)] If for any $L\in\mathcal{F}$ such that $L\subset \mathrm{Dom}(F)$ and $n(L)=0$, there exists an 
        $\mathcal{F}$-saturated neighborhood $S$ of $L$ in $\Sigma$ such that the restricted map $F\vert_S$ is $C^1$.
	\item[(H2)] If for any $L_*\in\mathcal{F}$ such that $L_*\subset\mathrm{Dom}(F)$, $1\leq n(L_*)\leq\infty$ and
        $F^{n(L_*)}(L_*)\subset\mathrm{Dom}(F)$,
        there is a connected neighborhood $S\subset\mathrm{Dom}(F)$ of $L_*$ such that the connected components
        $S_1$ and $S_2$ of $S\backslash L_*$ (possibly equal if $L_*\subset L_-\cup L_+$) satisfy the properties below:
        \begin{enumerate}
 		\item Both $F(S_1)$ and $F(S_2)$ are contained in $\Sigma\backslash (L_-\cup L_+)$.
 		\item For all $j\in\{ 1,2\}$, there exists a number $n^j(L_*)$ such that $1\leq n^j(L_*)\leq n(L_*)+1$
                and if $y_l\in S_j$ is a sequence converging to $y\in L_*$, then $\{ F(y_l)\} $ is a sequence converging to 
                $F^{n^j(L_*)}(y)$.
                If $n^j(L_*)=1$, then $F$ is $C^1$ in $S_j\cup L_*$.
 		\item If $L_*\subset \Sigma\backslash (L_-\cup L_+)$ (and so $S_1\neq S_2$), then either
                $n^1(L_*)=1$ and $n^2(L_*)>1$ or $n^1(L_*)>1$ and $n^2(L_*)=1$.
	\end{enumerate}
\end{description}

\begin{figure}[h]
\begin{center}
\includegraphics*[width=15cm]{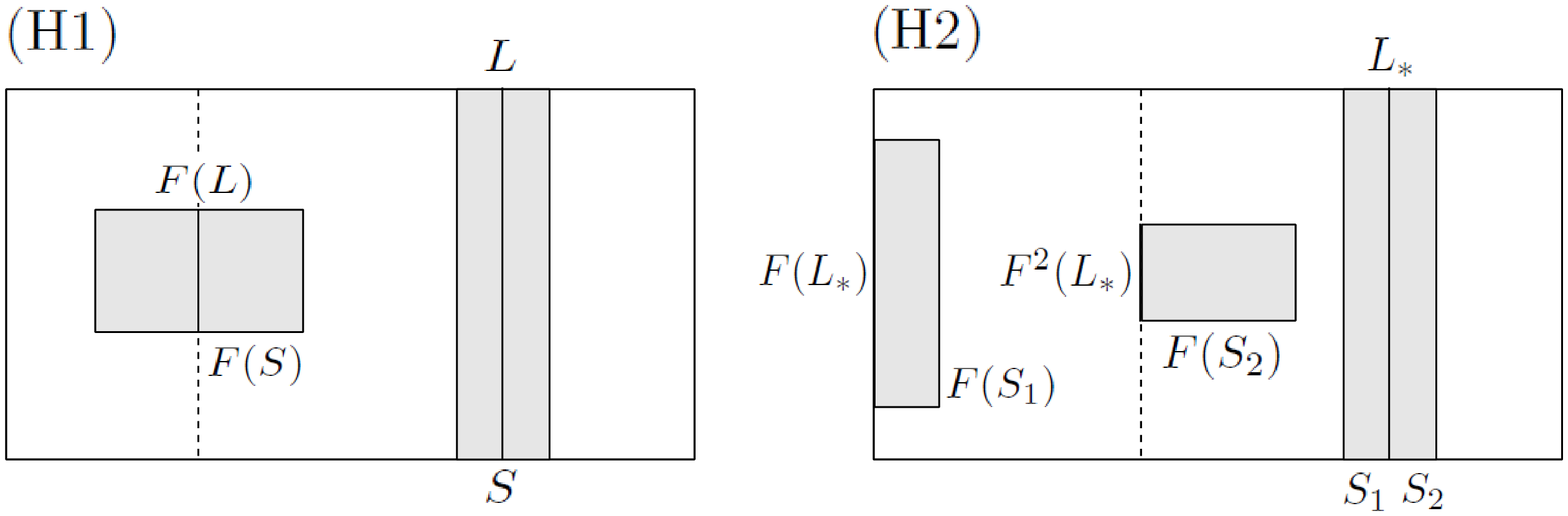}
\end{center}
\caption{(H1) and (H2) with $n(L_*)=1$, $n^1(L_*)=1$ and $n^2(L_*)=2$}
\end{figure}

\subsection{Proposition 1}
	This proposition gives a sufficient condition for the existence of periodic points of hyperbolic triangular maps.
	This proposition plays an important role in the proof of the Theorem.
\\

\noindent
\textbf{Proposition 1.}
Let $F$ be a $\lambda$-hyperbolic triangular map satisfying (H1), (H2), $\lambda>2$ and $\mathrm{Dom}(F)=\Sigma\backslash L_0$.
Then, $F$ has a hyperbolic periodic point.\\

	In [1,APPENDIX], the existence of periodic point was proved by a contradiction.
        Here, we prove that in a constructive way.
        Before the proof of Proposition 1, we need several preparations.
        Assume that $F$ is a $\lambda$-hyperbolic triangular map satisfing (H1), (H2), $\lambda>2$ and $\mathrm{Dom}(F)=\Sigma\backslash L_0$
        which contains the property: $L_-\cup L_+\subset\mathrm{Dom}(F)$.
        If there is a periodic point in $L_-\cup L_+$, then,
        it is hyperbolic because of a vertical contraction and $\lambda$-hyperbolicity of $F$.
        Therefore, we can suppose that there are no periodic points in $L_-\cup L_+$.
        Then, for a leaf $L\in\mathcal{F}$ with $L\subset\mathrm{Dom}(F)$, we have $0\leq n(L)\leq 2k$,
        for otherwise, there would exist a periodic point in $L_-\cup L_+$ because the iteration of $L$ passes $L_-\cup L_+$ at least 2k times.
        Let us call the following property the Hypotheses(*).\\
        \\
	Hypotheses(*): $F$ is a $\lambda$-hyperbolic triangular map satisfing (H1), (H2), $\lambda>2$,
        $\mathrm{Dom}(F)=\Sigma\backslash L_0$ and $\mathrm{Per}(F)\cap (L_-\cup L_+)=\emptyset$.\\
	\\
        We say that $F$ has the \textit{large domain} if $\mathrm{Dom}(F)=\Sigma\backslash L_0$.
\\

Let $k$ be the number of components of $\Sigma$.
We denote the leaf space of a vertial foliation $\mathcal{F}$ on $\Sigma$ by $SL$.
It is a disjoint union of $k$ copies of $I_1,\ldots ,I_k$ of $I$.
	For the 1-dimensional map $f:\mathrm{Dom}(f)\subset SL\longrightarrow SL$ induced by $F$,
        we define $f(L_*+):=\lim_{L\rightarrow L_*+}f(L)$ and $f(L_*-):=\lim_{L\rightarrow L_*-}f(L)$.
\\

\noindent
\textbf{Definition 10.}
	For $f$ induced by $F$ we define

\begin{enumerate}
	 \item $\mathcal{V}=\{ f(B); \; B\in\mathcal{F},\: B\subset\mathrm{Dom}(F)\cap(L_-\cup L_+)\} $
	 \item $\mathcal{L}_-=\bigcup\{ f(L_{0i}-); \; i\in \{ 1,\ldots ,k\} $ for which $f(L_{0i}-)$ exists.$\} $
	 \item $\mathcal{L}_+=\bigcup\{ f(L_{0i}+); \; i\in \{ 1,\ldots ,k\} $ for which $f(L_{0i}+)$ exists.$\} $
\end{enumerate}

	Define the \textit{discontinuous set of F} as
        $D(F)=\{ x\in\mathrm{Dom}(F);\; F$ is discontinuous at $x\}$.
We denote the foliation and the cone field associated to $F$ by $\mathcal{F}$ and $C_{\alpha}$ respectively.
Let $<$ be the natural order in the leaf space $I_i$ of $\mathcal{F}_i$, where $\mathcal{F}_i$ is a 
vertical foliation in $\Sigma_i$ ($i=1,\ldots ,k$).
A \textit{vertical band} in $\Sigma$ is a region between two disjoint vertical curves $L$ and $L'$ in the
same component $\Sigma_i$ of $\Sigma$.
The notation $[L,L']$ and $(L,L')$ indicates closed and open vertical band respectively.
	If $c$ is a curve in $\Sigma$, we denote its end points by $c_0,c_1$,
        its closure $Cl(c)=c\cup\{ c_0,c_1\} $, and its interior $\mathrm{Int}(c)=c\backslash \{c_0,c_1\} $.
An open curve is a curve without end points.
We say that $c$ is \textit{tangent} to $C_{\alpha}$ if $c'(t)\in C_{\alpha}(c(t))$ for all $t\in\mathrm{Dom}(c)$.
	The next lemma is proved in \cite{BM}.
\\

\noindent
\textbf{Lemma 1} [1,Lemma 14].
For every open curve $c\subset \mathrm{Dom}(F)\backslash D(F)$ tangent to $C_{\alpha}$ there are an open curve
$c^*\subset c$ and $n'(c)>0$ such that $F^j(C^*)\subset \mathrm{Dom}(F)\backslash D(F)$ whenever $0\leq j\leq n'(c)-1$ and 
$F^{n'(c)}(c^*)$ covers a band $(W,W')$ with
	\[ W,W'\subset L_- \cup L_+ \cup \mathcal{V} \cup \mathcal{L}_- \cup \mathcal{L}_+. \]

	Now let us prove proposition 1.
\\

\noindent
\textit{Proof of Proposition 1.}
        As we mensioned before, if there is a periodic point in $L_-\cup L_+$, Proposition 1 is proved.
        It remains to prove is the existence of a hyperbolic periodic point under the Hypothesis(*).
	It is known that $\mathrm{Dom}(F)\backslash D(F)$ is open-dense in $\Sigma$ (See \cite{BM}.).
	Define
        \[	\mathcal{B} = \{ (W,W');\; W,W' \subset L_- \cup L_+ \cup \mathcal{V} \cup \mathcal{L}_- \cup \mathcal{L}_+\}. \]
It is clear that $\mathcal{B}$ is a finite set.
In $\mathcal{B}$ we define the relation $B \leq B'$ if and only if there are an open curve $c\subset B$ 
tangent to $C_{\alpha}$ with $\mathrm{Cl}(c)\subset \mathrm{Dom}(F)\backslash D(F)$, 
an open subcurve $c^*\subset c$ and $n>0$ such that
$F^j(c^*)\subset \mathrm{Dom}(F)\backslash D(F)$ whenever $0\leq j\leq n-1$
and $F^n(c^*)$ covers $B'$.

As $\mathrm{Dom}(F)\backslash D(F)$ is open-dense in $\Sigma$ and the bands in $\mathcal{B}$ are open, 
we can use Lemma 1 to prove that for every $B\in \mathcal{B}$ there is $B'\in \mathcal{B}$ such that $B \leq B'$.
Then, we can construct a chain
\[ B_{j_1} \leq B_{j_2} \leq B_{j_3} \leq \ldots \; (j_i\in\{ 1,\ldots ,m\}) \]
As $\mathcal{B}$ is finite it would exist a closed sub-chain
\[ B_{j_i} \leq B_{j_{i+1}} \leq \ldots \leq B_{j_{i+s}} \leq B_{j_i} \]
Hence there is a positive integer $n$ such that $F^n(B_{j_i})$ covers $B_{j_i}$.
	Since $F$ preserves $\mathcal{F}$, there exists a leaf $\tilde{L}\in\mathcal{F}$ such that
        $F^n(\tilde{L})\subset\tilde{L}$, implying the existence of a periodic point in $\tilde{L}$.
        By a vertical contruction and a horizontal expansion which is derived by the $\lambda$-hyperbolicity 
        of $F$, this periodic point is hyperbolic, and therefore Proposition 1 is proved.
\wsq

\subsection{Singular-cross section and induced foliation}
	In this subsection we construct a family of cross-sections and foliations on them.
Let $X_t$ be a $C^1$ flow and let $\sigma$ be a Lorenz-like singularity of $X_t$.
Then, $\sigma$ is hyperbolic,
	and we have invariant manifolds $W_X^s(\sigma)$, $W_X^u(\sigma)$ and $W_X^{ss}(\sigma)$ with
	$\mathrm{dim}(W_X^s(\sigma))=2$, $\mathrm{dim}(W_X^u(\sigma))=1$ and $\mathrm{dim}(W_X^{ss}(\sigma))=1$.
        (See \cite{BM}.)

$W_X^{ss}(\sigma)$ separates $W_X^s(\sigma)$ into two connected components, namely, the top one and the bottom one.
In the top component, we consider a cross-section $S_{\sigma}^t$ of $X_t$ together with a curve $l_{\sigma}^t$
	which goes directly to $\sigma$.
Similarly we consider a cross-section $S_{\sigma}^b$ and a curve $l_{\sigma}^b$ in the bottom component.
We take the section $S_{\sigma}^*$ to be diffeomorphic to $[-1,1]\times [-1,1]$ and the curve $l_{\sigma}^*$ to be 
contained in $W_X^s(\sigma)\backslash W_X^{ss}(\sigma)$ for $*=t,b$.
The positive orbits of $X_t$ starting at $(S_{\sigma}^t\cup S_{\sigma}^b)\backslash (l_{\sigma}^t\cup l_{\sigma}^b)$ 
exit a small neighborhood of $\sigma$ passing through the cusp region.
The positive orbits starting at $l_{\sigma}^t\cup l_{\sigma}^b$ goes directly to $\sigma$.
The boundary of $S_{\sigma}^*$ is formed by four curves, two of them transverse to $l_{\sigma}^*$ 
and two of them parallel to $l_{\sigma}^*$.
The union of the curves in the boundary of $S_{\sigma}^*$ which are parallel (resp. transverse) to $l_{\sigma}^*$ 
is denoted by $\partial^vS_{\sigma}^*$ (resp. $\partial^hS_{\sigma}^*$).
        The cross-sections $S_{\sigma}^t$ and $S_{\sigma}^b$ above are called \textit{singular cross-sections associated to $\sigma$}.
	The curves $l_{\sigma}^t$ and $l_{\sigma}^b$ are called \textit{singular curves} of $S_{\sigma}^t$ and $S_{\sigma}^b$ respectively.
        We also call a family of disjoint cross-sections
	$S=\{ S_{\sigma}^t,S_{\sigma}^b;\; \sigma\in \mathrm{LSing}_X(\Lambda)\} $ with $\Lambda\cap\partial^hS=\emptyset$ 
        the \textit{singular cross-section} of $\Lambda$.
        Similarly we call a family of disjoint curves
        $l=\{ l_{\sigma}^t,l_{\sigma}^b;\; \sigma\in \mathrm{LSing}_X(\Lambda)\} $
        the \textit{singular curves} of $S$.
        Define
        \[ \partial^hS=\bigcup_{\sigma\in\mathrm{LSing}_X(\Lambda)}(\partial^hS_{\sigma}^t\cup \partial^hS_{\sigma}^b) \quad \mathrm{and} \quad
        \partial^vS=\bigcup_{\sigma\in\mathrm{LSing}_X(\Lambda)}(\partial^vS_{\sigma}^t\cup \partial^vS_{\sigma}^b). \]

\begin{figure}[h]
\begin{center}
\includegraphics*[width=12cm]{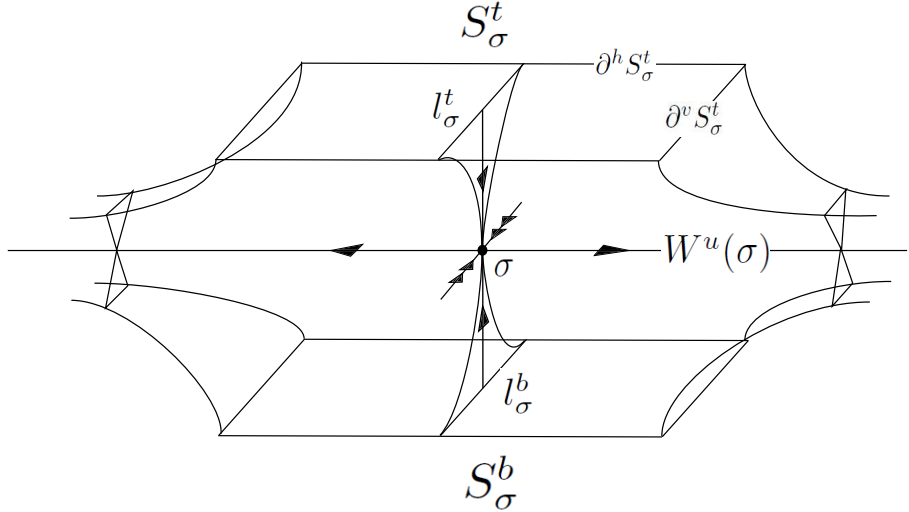}
\end{center}
\caption{The singular cross-section and the singular curve}
\end{figure}

Now we construct a foliation on the singular cross-section.
For the singular-hyperbolic splitting
$T_{\Lambda}M=E_{\Lambda}^s\oplus E_{\Lambda}^c$,
$E^s_{\Lambda}$ and $E^c_{\Lambda}$ can be extended continuously to invariant splittings
$E^s_{U(\Lambda)}$ and $E^c_{U(\Lambda)}$ on a neighborhood $U(\Lambda)$ of $\Lambda$, respectively.
In particular, the contracting direction is 1-dimensional.
The standard Invariant Manifold Theorem\cite{HPS} implies that $E_{U(\Lambda)}^s$ is integrable,
i.e. tangent to an invariant continuous 1-dimensional contracting foliation $\mathcal{F}^{ss}$ on $U(\Lambda)$.
Let $S$ be a singular cross-section of $\Lambda$ contained in $U(\Lambda)$.
	We construct a foliation $\mathcal{F}$ on $S$ by projecting $\mathcal{F}^{ss}$ onto $S$ along the flow.
	(See \cite{BM} for the precise construction.)

%%------------------------------------------------------------------------------------------------------------------
%%------------------------------------------------------------------------------------------------------------------

\section{PROOF OF THE THEOREM}
	In this section, we prove the Theorem.
        For a singular-hyperbolic Lyapunov stable set $\Lambda$,
        first we consider two exceptional cases where there are no singularities in $\Lambda$, and
        the case where there are singularities except for Lorenz-like ones.

\subsection{The exceptional cases}

\noindent
	\textbf{Proposition 2.}
	Let $\Lambda$ be a singular-hyperbolic Lyapunov stable set.
	If there are no singularities in $\Lambda$, then $\Lambda$ has a hyperbolic periodic orbit.
\\

\noindent
\textit{Proof.}
	$\Lambda$ is a (saddle-type) hyperbolic set.
	Take $x\in\Lambda$ and a cross-section $\Sigma$ with $x$.
	By the Shadowing Lemma\cite{HK}, there exists a periodic point $p$ near $x$ and hence near $\Lambda$.
	Assume $p\notin\Lambda$ and take $U\supset\Lambda$ such that $p\notin U$.
	The stable and the unstable manifolds of $x$ and $p$ are large enough to intersect transversally.
	By the $\lambda$-lemma\cite{PM}, any neighborhood of $x\in V_{Lyp}(U)$ can be arbitarily close to $p$ under some iteration,
	contradicting the Lyapunov stability (Figure 3).
	Therefore we obtain a periodic point in $\Lambda$, moreover this is hyperbolic because of a 
        contraction and an expansion derived from the singular-hyperbolicity of $\Lambda$.
\wsq\\

\begin{figure}[h]
\begin{center}
\includegraphics*[width=6cm]{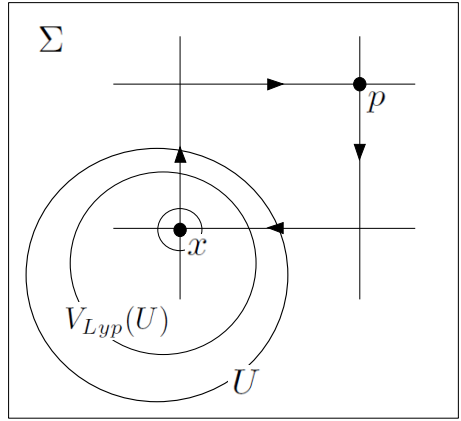}
\end{center}
\caption{The Lyapunov stability and the transversal intersection of invariant manifolds.}
\end{figure}

\noindent
	\textbf{Proposition 3.}
	Let $\Lambda$ be a singular-hyperbolic Lyapunov stable set of a $C^1$ flow $X_t$.
	If $\Lambda$ has singularities except for Lorenz-like ones, then $\Lambda$ has a hyperbolic periodic orbit.
\\

\noindent
\textit{Proof.}
	We can take $x\in\Lambda\backslash\mathrm{Sing}(X)$.
        For otherwise, $\Lambda$ would be a set of singularities and they are all hyperbolic which are discrete.
        Each singularity has both positive and negative eigenvalues because $\Lambda$ is singular-hyperbolic,
        contradicting the fact that $\Lambda$ is Lyapunov stable.
	Clearly $\omega_X(x)\subset\Lambda$ since $\Lambda$ is compact invariant.
	Let us see that $\omega_X(x)$ has no singularities.
        If there exists $\sigma\in\mathrm{Sing}_X(\omega_X(x))$,
        it is a singularity of Lorenz-like or one with two positive eigenvalues.
        Here, we have assumed there are no Lorenz-like singularities, it has two positive eigenvalues.
        (For such singularity, $W^{ss}_X(\sigma)=W^s_X(\sigma)$.)
        Then, there exists $q\in\omega_X(x)\cap(W^s_X(\sigma)\backslash\{\sigma\})$ (Figure 4).
        
        \begin{figure}[h]
	\begin{center}
	\includegraphics*[width=9cm]{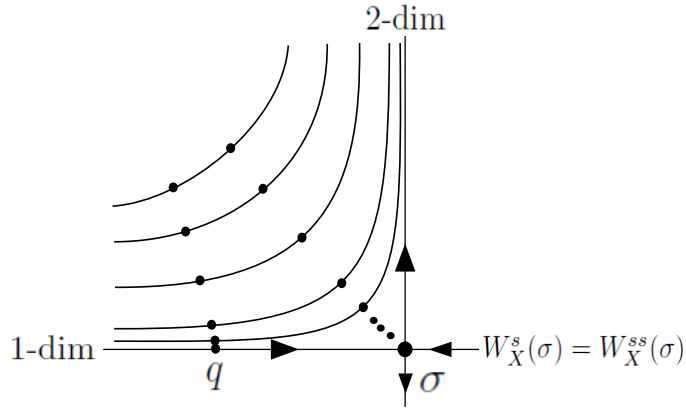}
	\end{center}
        \caption{The case where $\omega_X(x)$ contains a singularity with two positive eigenvalues.}
	\end{figure}
        
        This and $\omega_X(x)\subset\Lambda$ contradict the fact that $\Lambda\cap(W^s_X(\sigma)\backslash\{\sigma\})=\emptyset$.
        Since $\omega_X(x)\subset\Lambda$ and $\omega_X(x)$ has no singularities, $\omega_X(x)$ is a (saddle-type) hyperbolic set.
	Then, by the same argument as in the proof of Proposition 2,
        we obtain a hyperbolic periodic point in $\Lambda$.
\wsq\\

	Now we consider the case where there are Lorenz-like singularities in $\Lambda$.

\subsection{Preliminaries for the proof.}
For the proof, we need a lemma dealing with the return maps associated to singular cross-sections.
Let $\Lambda$ be a singular-hyperbolic Lyapunov stable set of a $C^1$ flow $X_t$.
Associated to any singular cross-section $S$ of $\Lambda$, we have a return map
$\Pi=\Pi_S:\mathrm{Dom}(\Pi)\subset S\rightarrow S$
given by $\Pi(x)=X_{T(x)}(x)$ where $T(x)$ denotes the first positive return time of $x$.\\

Here, using the foliation $\mathcal{F}$ of subsection 2.4,
we \textit{refine} a singular cross-section $S\subset U(\Lambda)$ in the following way.
Let $S$ be a singular cross-section of $\Lambda$.
By the construction, $l_{\sigma}^*$ divides $S_{\sigma}^*$ into two connected components $S_{\sigma}^{*,+}$ and $S_{\sigma}^{*,-}$ ($*=t,b$).
For a small $\delta>0$, we choose two points $x_{\delta}^+,x_{\delta}^- \in S_{\sigma}^{*,\pm}$ whose distance
to $l_{\sigma}^*$ is $\delta$.
Define $S_{\sigma}^*(\delta)$ as the singular cross-sections of $\sigma$ satisfying the following property:
	\[ \partial^vS_{\sigma}^*(\delta)=\mathcal{F}_{x_{\delta}^-}\cup \mathcal{F}_{x_{\delta}^+}. \]
Since $S$ is a singular cross-section of $\Lambda$, we conclude that the set
	\[ S(\delta)=\{ S_{\sigma}^t(\delta),S_{\sigma}^b(\delta);\; \sigma\in \mathrm{LSing}_X(\Lambda)\} \]
is also a singular cross-section of $\Lambda$.
Note that $S$ and $S(\delta)$ have the same singular curve $l$.

We also refine its return map.
For the refinement $S(\delta)$, we denote the return map associated to $S(\delta)$ by $\Pi_{\delta}=\Pi_{S(\delta)}$
and denote  the return time of $x\in \mathrm{Dom}(\Pi_{\delta})$ by $T_{\delta}(x)$.
Clearly $S(\delta)\subset S$ and so $S(\delta)\subset U(\Lambda)$ for all $\delta$.
A simple but important observation is that the return time $T_{\delta}$ is uniformly large as $\delta\rightarrow 0^+$,
	\[ \lim_{\delta\rightarrow 0^+}\inf_{x\in S(\delta)}T_{\delta}(x)=\infty . \]

For a singular cross-section $S$ of $\Lambda$ and its refinement $S(\delta)$,
note that each components of $S(\delta)$ can be identified with the square $I^2=I\times I$ (where $I=[-1,1]$) 
such that its singular curve corresponds to $\{ 0\} \times [-1,1]$ and its vertical boundaries correspond to 
$\{\pm 1\} \times [-1,1]$.
It follows that $S(\delta)$ can be identified with a finite collection $\Sigma=\Sigma_{\delta}$
of squares with
	\[ \partial^vS(\delta)=L_-\cup L_+ \quad \mathrm{and} \quad l=L_0 \]
where $\Sigma$, $L_-$, $L_+$ and $L_0$ are as in subsection 2.1.

With these identifications, we define
	\[ F=\Pi_{\delta} \quad \mathrm{and} \quad \mathrm{Dom}(F)=\mathrm{Dom}(\Pi_{\delta}).\]
Of course $F$ and Dom($F$) depend on $\delta$.
Hence we have a map
	\[ F:\mathrm{Dom}(F)\subset \Sigma\rightarrow \Sigma \]
which is the return map induced by the flow $X_t$ on the section $S(\delta)$.\\

It is clear that $\mathrm{Dom}(\Pi)\subset S\backslash l$
where $l$ is the singular curve of $S$.
	We say that $\Pi$ has the \textit{large domain} if $\mathrm{Dom}(\Pi)=S\backslash l$.
	The following lemma proves that there exists a singular cross-section whose return map has the large domain and satisfies 
        some conditions simultaneously, if the flow $X_t$ has no periodic orbits in $\Lambda$.
\\

\noindent
\textbf{Lemma 2.}
Let $\Lambda$ be a singular-hyperbolic Lyapunov stable set of $X_t$ and let $\lambda>0$ be fixed.
If $X_t$ has no periodic orbits in $\Lambda$, then, for any neighborhood $U\supset\Lambda$,
there exists a singular cross-section $S\subset U$ such that
its return map $\Pi$ is a $\lambda$-hyperbolic triangular map
satisfying (H1), (H2) and that $\mathrm{Dom}(\Pi)=S\backslash l$.\\

\noindent
\textit{Proof.}
	First, define
        \[ \mathrm{Sing}^*=\mathrm{Sing}_X(\Lambda)\backslash\mathrm{LSing}_X(\Lambda) \quad \mathrm{and} \quad
        W=\bigcup_{\sigma\in\mathrm{Sing}_X(\Lambda)}W^{ss}_X(\sigma). \]
        Since $W^{ss}_X(\sigma)\cap\Lambda=W^s_X(\sigma)\cap\Lambda=\{ \sigma\}$ for $\sigma\in\mathrm{Sing}^*$ and
        $W^{ss}_X(\sigma)\cap\Lambda=\{ \sigma\}$ for $\sigma\in\mathrm{LSing}_X(\Lambda)$,
        we have that $\Lambda\cap(W\backslash\mathrm{Sing}_X(\Lambda))=\emptyset$.
        
        Let us exhibit a contradiction assuming that there exists $U\supset\Lambda$ such that
        any singular cross-section $S\subset U$ whose return map is a $\lambda$-hyperbolic triangular map satisfying (H1) and (H2)
        has a point $x\in S\backslash l$ satisfying $x\notin\mathrm{Dom}(\Pi)$.
	Here we can take such $U$ with $U\subset U(\Lambda)$ and $W\not\subset U$.
        Moreover we can assume $U\cap(\mathrm{Sing}(X)\backslash\mathrm{Sing}_X(\Lambda))=\emptyset$.
        For otherwise, $\mathrm{Sing}(X)\backslash\mathrm{Sing}_X(\Lambda)$ accumulates on $\Lambda$,
        and letting $\tilde{\sigma}$ be an accumulation point, we have $\tilde{\sigma}\in\mathrm{Sing}_X(\Lambda)$,
        which is hyperbolic by the definition of singular-hyperbolicity.
        However, this contradicts Grobman-Hartman Theorem\cite{PM}.
        
        The following property is known.\\
        
        \noindent
	[1, Lemma 4 and Proposition 2]
        For any $U\supset\Lambda$, there is a singular cross-section $S\subset U$ associated to $\Lambda$
        which has a small diameter and is close to $\Lambda$ such that
	if $S(\delta)$ is the refinement of $S$, then for all small $\delta>0$,
	$F=\Pi_{\delta}$ is a $\lambda$-hyperbolic triangular map satisfying (H1) and (H2).\\

        Let us take a sequence $\{S_k\}$ of such singular cross-sections accumulating on $\Lambda$.
        For each $S_k$, take $\delta_k>0$ so that $S_k(\delta_k)$ is a refinement of $S_k$.
        Then, we consider the sequence $\{S_k(\delta_k)\}$ of refined singular cross-sections satisfying the following property:
        for Lorenz-like singularities
        (Here for simplicity, in their neighborhoods, we identify eigenspaces of Lorenz-like singularities with (x,y,z)-axes, respectively, as depicted in Figure 5)
        $S_k(\delta_k)$ contains the rectangular region $C_k$ in $\{z=d_k\}$ satisfying that
        $\mathrm{dist}(\partial^hC_k,(0,0,d_k))=d_k$ and $\mathrm{dist}(\partial^vC_k, l_k)=\epsilon_k<\delta_k$
        (where $d_k=\mathrm{dist}(\mathrm{LSing}_X(\Lambda),S_k(\delta_k))$ and $l_k$ is the singular curve of $S_k(\delta_k)$),
        and if
        \[	\{\gamma_i\}=\left\{ \bigcup_{0\leq t\leq t_i}X_t(x_i),\; x_i\in V_{Lyp}(U) \right\}	\]
        accumulates on a Lorenz-like singularities, then $\gamma_i\cap S_k(\delta_k)\in C_k$ for every large $i\in\NNN$.
        
        \begin{figure}[h]
        \begin{center}
	\includegraphics*[width=6cm]{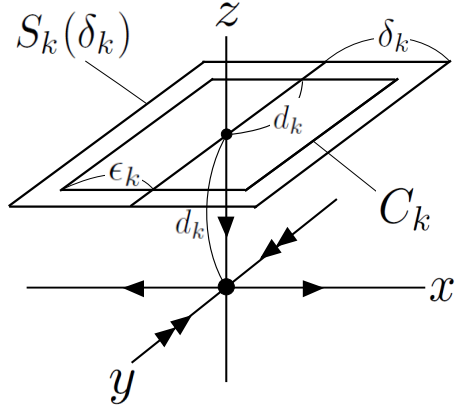}
	\caption{Singular cross-section $S_k(\delta_k)$ and $C_k$.}
	\end{center}
        \end{figure}
        
        Now, we take a family $\{U_n\}$ of neighborhoods of $\Lambda$ satisfying $U\supset U_0\supset U_1\supset \cdots \supset \Lambda$.
	For each $U_n$, take $V_{Lyp}(U_n)$ and a singular cross-section $S_{k_n}(\delta_{k_n})\subset V_{Lyp}(U_n)$ from the above sequence.
        Let $l_n$ be a singular curve of $S_{k_n}(\delta_{k_n})$ and $\Pi_n$ the return map of $S_{k_n}(\delta_{k_n})$.
        By the hypothesis, there exists $x_n\in S_{k_n}(\delta_{k_n})\backslash l_n$ such that $x_n\notin\mathrm{Dom}(\Pi_n)$.
        We note that $\omega_X(x_n)\cap(\mathrm{Sing}(X)\backslash\mathrm{Sing}_X(\Lambda))=\emptyset$
        since $U\cap(\mathrm{Sing}(X)\backslash\mathrm{Sing}_X(\Lambda))=\emptyset$ and $\omega_X(x_n)\subset \overline{U}_n\subset U.$
        (Here $\overline{U}_n$ is the closure of $U_n$.)
        
        Let us see $\omega_X(x_n)\cap\mathrm{Sing}_X(\Lambda)=\emptyset$.
        If it is not so, $\omega_X(x_n)$ contains a singularity of $\mathrm{Sing}^*$ or $\mathrm{LSing}_X(\Lambda)$.
        In the case where $\exists\sigma\in\omega_X(x_n)\cap\mathrm{Sing}^*$,
        there exists $q\in\omega_X(x_n)\cap (W^s_X(\sigma)\backslash\{\sigma\})$ with $q\notin U$, contradicting the Lyapunov stability.
        In the case where $\exists\sigma\in\omega_X(x_n)\cap\mathrm{LSing}_X(\Lambda)$,
        we have $X_T(x_n)\in C_{k_n}\subset S_{k_n}(\delta_{k_n})$ for a large $T>0$ by the construction of $S_{k_n}(\delta_{k_n})$.
        This contradicts that $x_n\notin\mathrm{Dom}(\Pi_n)$.
        
        For $\{\omega_X(x_n)\}$, let $P$ be the set of accumulation points of $\{\omega_X(x_n)\}$, then $P\subset\Lambda$.
        We have that $\mathrm{Sing}_X(P)\neq\emptyset$, for otherwise $P$ would be hyperbolic and the
        same argument as in Proposition 2 would lead a contradiction.
        Then, we have two cases as before.
        In the case where $\exists\sigma\in P\cap\mathrm{Sing}^*$,
        there exists $q\in P\cap (W^s_X(\sigma)\backslash\{\sigma\})$ with $q\in\Lambda$ since $P\subset\Lambda$,
        which contradicts that $\Lambda\cap (W\backslash\mathrm{Sing}_X(\Lambda))=\emptyset$.
        In the case where $\exists\sigma\in P\cap\mathrm{LSing}_X(\Lambda)$,
        we have $\omega_X(x_N)\cap S_{k_N}(\delta_{k_N})\in C_{k_N}$ for a large $N\in\NNN$,
        then there exists $T>0$ such that $X_T(x_N)\in S_{k_N}(\delta_{k_N})$.
        This contradicts $x_N\notin\mathrm{Dom}(\Pi_N)$.
        
        Thus, Lemma 2 has been proved.\wsq

\subsection{Proof of the Theorem.}	
        Two exceptional cases have been proved before; that is the case where $\Lambda$ contains no singularities
        and the case where $\Lambda$ contains singularities with two positive eigenvalues but no Lorenz-like ones.
        Therefore we can suppose that $\Lambda$ contains at least one Lorenz-like singularity.
	Let us assume that there are no periodic orbits in $\Lambda$.
        By Lemma 2, there exists a singular cross-section close to $\Lambda$ such that
        the return map $F$ is a $\lambda$-hyperbolic triangular map satisfying (H1) and (H2) moreover $F$ has the large domain.
        Then $F$ satisfies the hypothesis of Proposition 1.

        Let $\{ c_n\} $ be a family of curves accumulating on $\Lambda$.
        Define
        \[ 	\mathcal{B}=\{ (W,W');\; W,W'\subset L_+ \cup L_- \cup \mathcal{V} \cup \mathcal{L}_+ \cup \mathcal{L}_- \} .	\]
	Then this is a finite set.
        As in the proof of Proposition 1, using Lemma 1, 
        we see that some $F$-iteration of $c_n$ covers $B_{n_1}\in\mathcal{B}$.
        Again by Lemma 1, some $F$-iteration of $B_{n_1}$ covers $B_{n_2}$.
        Repeating this process, we obtain a chain of elements of $\mathcal{B}$.
        Since $\mathcal{B}$ is a finite set, we obtain a closed sub-chain, which is called a cycle $\beta_{n}$.
        In this way, each $c_n$ coresponds to a cycle $\beta_{n}$.
	Since the set of cycles are also finite, there exists a cycle $\bar{\beta}$ that is an accumulation cycle of $\{ \beta_n\}$.
        Let $\{\tilde{c}_n \} $ be a subsequence of $\{ c_n\}$ such that $\tilde{c}_n$ corresponds to $\bar{\beta}$.
        Let $p$ be a periodic point in the cycle $\bar{\beta}$.
        Assume that $p\notin\Lambda$ and take a neighborhood $U\supset\Lambda$ with $p\notin U$.
        For $\{ \tilde{c}_n \} $, take an integer $N$ large enough to satisfy $\tilde{c}_N\subset V_{Lyp}(U)$.
        Then, some iteration of $\tilde{c}_N$ under $F$ covers $p$.
        Note that the size of the stable manifold of $p$ is bigger than the length of the section because
        we have considered a foliation on the section in $U(\Lambda)$ introduced in subsection 2.4.
        Take an integer $m$ for which $F^m(\tilde{c}_N)$ covers $p$.
        Since $F$ is $\lambda$-hyperbolic, $F^m(\tilde{c}_N)$ has a transversal intersection with $W^s(p)$.
        Let $q$ be the intersection point of them, then, the iteration images of $q$ under $F$ accumulates on $p$.
	Since $p\notin U$, this implies that some iteration of $q\in F^m(\tilde{c}_N)$ and therefore that of $F^{-m}(q)\in \tilde{c}_N\subset V_{Lyp}(U)$
        is not contained in $U$.  
        This contradicts the Lyapunov stability, proving that $p\in\Lambda$.
\wsq

%%------------------------------------------------------------------------------------------------------
%%------------------------------------------------------------------------------------------------------

\section*{APPENDIX}
We give an example of a Lyapunov stable set which is not attracting by modifing the GLA and using the Cherry-flowbox\cite{PM}.
This example is not singular-hyperbolic, however its property is close to singular-hyperbolicity.

The Cherry-flow is a vector field on the 2-torus with one sink and one saddle. Both singularities are hyperbolic.
We put a saddle $s$, and a sink $p$.
See [6,APPENDIX] for the construction and properties of the Cherry-flow.
Identifying the 2-torus with $[0,1] \times [0,1] $, we depict the Cherry-flow in Figure 6 (left).
The time-reversed flow is depicted in Figure 6 (center).
For the time-reversed Cherry-flow, we assume eigenvalues $\lambda_-$ and $\lambda_+$ of the saddle $s$ satisfying
$\lambda_-<0<-\lambda_-=\lambda_+$.

\begin{figure}[h]
\includegraphics*[width=5cm]{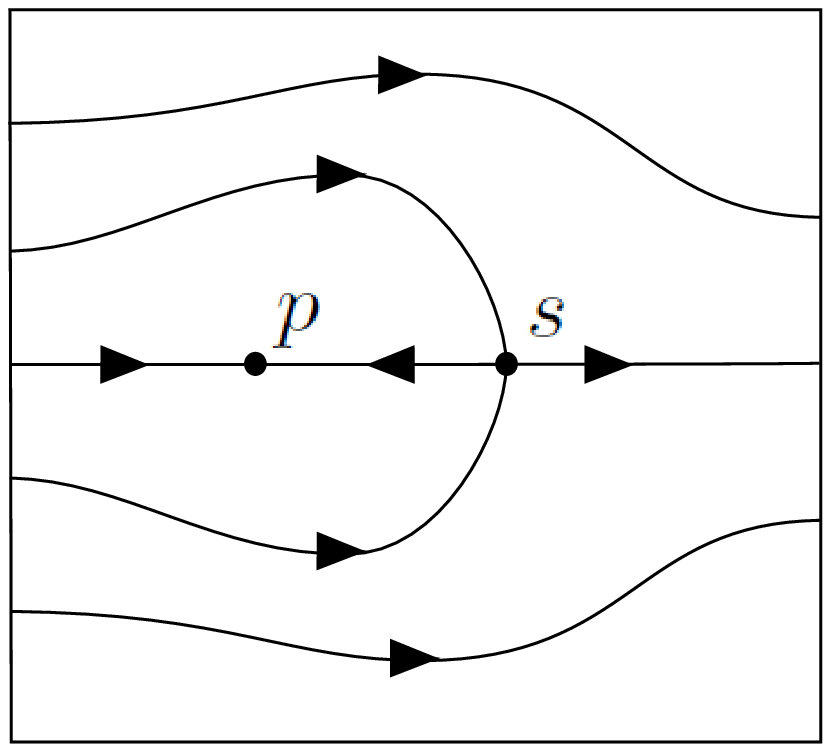}
\includegraphics*[width=5cm]{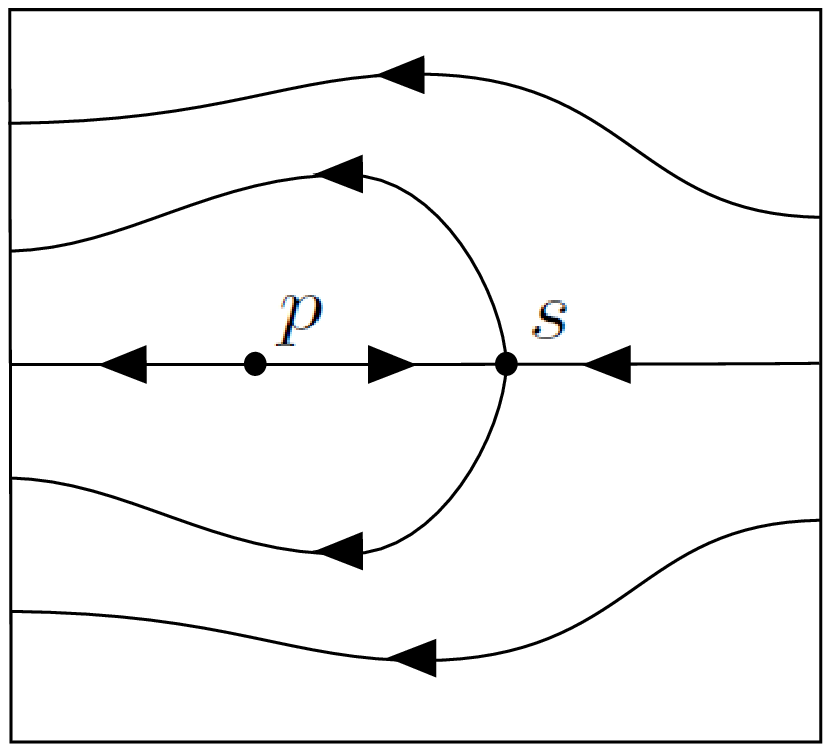}
\includegraphics*[width=5cm]{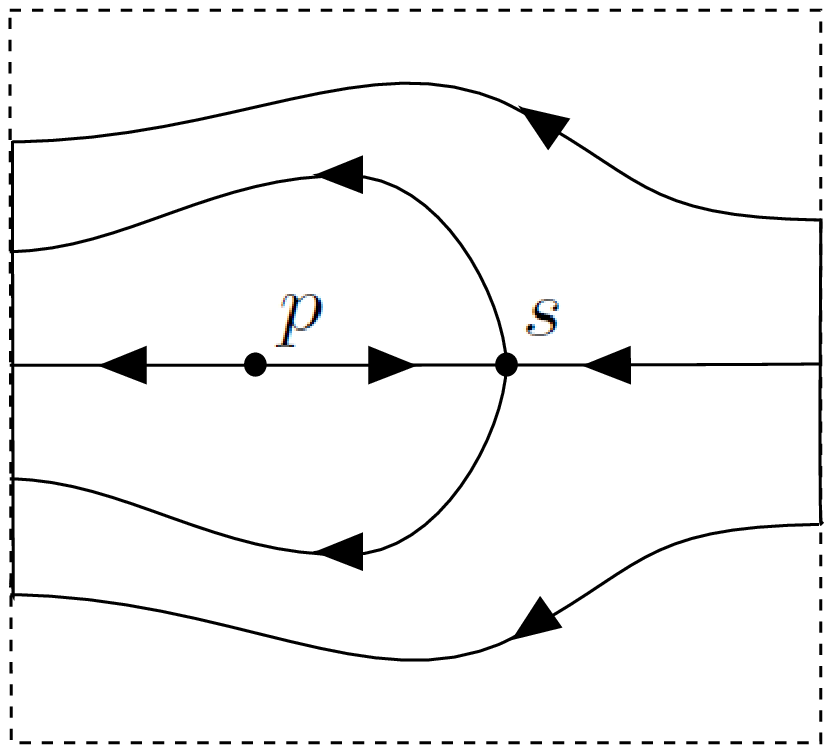}
\caption{The Cherry-flow and its arrangements}
\end{figure}

Taking a part of the center figure (Figure 6, right) and
multipling a contracting direction, we obtain the Cherry-flowbox $C$ (Figure 7).
Identify $C$ with $([0,1]\times[0,1])\times[-1,1]$.
Then, $C$ is devided into two components by $[0,1]\times\{\frac{1}{2}\}\times[-1,1]$ which we call $C^t$ and $C^b$.\\

\begin{figure}[h]
\begin{center}
\includegraphics*[width=10cm]{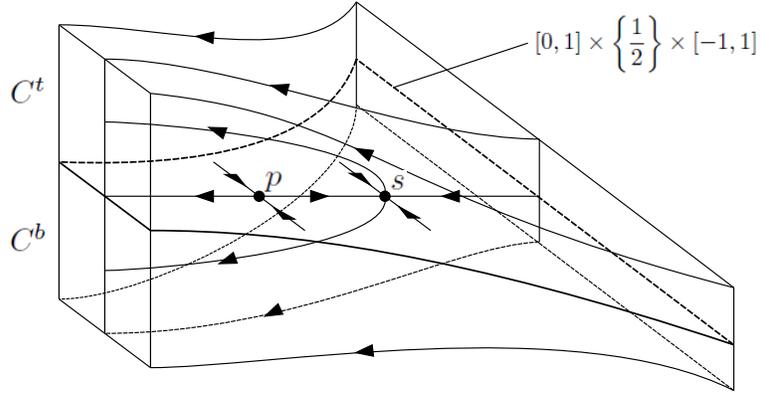}
\end{center}
\caption{Cherry-flowbox C}
\end{figure}

Now, let $\Sigma$ be a singular cross-section associated to a singularity $\sigma$ which has real eigenvalues
$\lambda_1$, $\lambda_2$ and $\lambda_3$ satisfying $\lambda_2<\lambda_3<0<-\lambda_3=\lambda_1$.
We call its two components $\Sigma^t$ and $\Sigma^b$.
We identify each of them with $[-1,1] \times [-1,1]$.
Moreover $\Sigma$ is devided into four regions by singular curves $l^t$ and $l^b$.
We put them as $\Sigma^t_-$, $\Sigma^t_+$, $\Sigma^b_-$ and $\Sigma^b_+$.
Let $\{ -1 \} \times \{ -1/3 \} \in \Sigma^t$ be the position of the first intersection with $\Sigma$
of one component of $W^u(\sigma)$, namely $W^u_+(\sigma)$.

Then, we take a flowbox around a part of the other component of $W^u(\sigma)$, namely $W^u_-(\sigma)$,
and replace it by the Cherry-flowbox $C$ with the contracting direction.
Also we connect $W^u_-(\sigma)$ with $W^s(s)$.
Then, the orbit of $W^u_-(\sigma)$ goes into $C$ and converges to the hyperbolic saddle $s$ in $C$.
On the other hand, two components of $W^u(s)$ go out of $C$.
Here we make an important remark.
Parts of the cusp regions coming from $\Sigma^t_-$ and $\Sigma^b_-$ enter into $C^t$ and $C^b$ respectively.
By the choices of $C^t$ and $C^b$, two cusp regions run along $W^u_-(\sigma)$ and $W^u(s)$ without intersecting each other and go out of $C$.
Now we set the returns of them.
We put the first intersecting points of $W^u_t(s):=W^u(s)\subset C^t$ and $W^u_b(s):=W^u(s)\subset C^b$ as
$\{1\}\times\{1/3\}\in \Sigma^t$ and $\{0\}\times\{0\}\in \Sigma^b$, respectively.

Using the Cherry-flowbox, we have constructed a vector field depicted in Figure 8.
Then, define
\[ \Lambda = \left( \bigcap_{T\geq 0}\bigcup_{t\geq T}X_t(\Sigma^t)\right) \cup W^u(\sigma) \cup W^u(s) \cup \{ \sigma,s \}. \]

\begin{figure}[h]
\begin{center}
\includegraphics*[width=10cm]{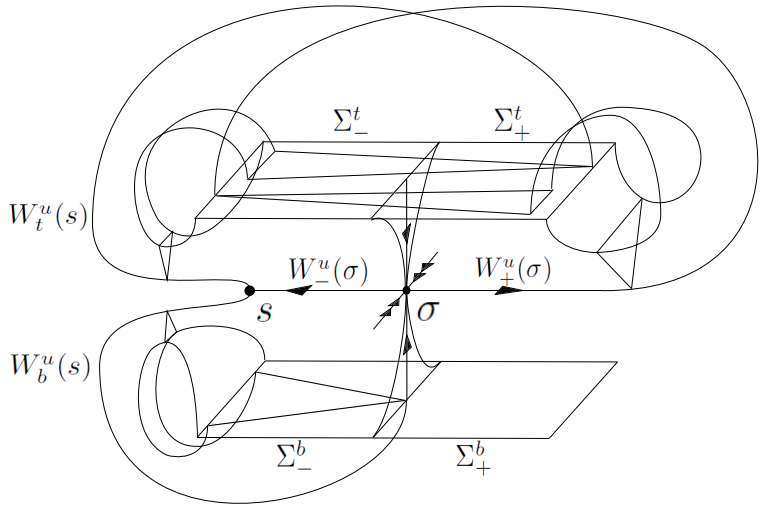}
\end{center}
\caption{vector field}
\end{figure}

The return map $F$ on $\Sigma\backslash (l^t\cup l^b)$ consists of the following two parts:
$F^t:\Sigma^t\backslash l^t\longrightarrow\Sigma^t$ and 
$F^b_-:\Sigma^b_-\backslash l^b\longrightarrow\Sigma^b_-$.
We can assume that the return map $F$ is a triangular-map.
So, $F$ is reduced to 1-dimensional map $f$.
As we mensioned before, both $\Sigma^t$ and $\Sigma^b$ are identified with $[-1,1]\times[-1,1]$,
and 1-dimensional maps to which $F^t$ and $F^b_-$ are reduced are
$f^t:[-1,1]\longrightarrow[-1,1]$ and
$f^b_-:[-1,0]\longrightarrow[-1,0]$, respectively.
Here, we assume the following conditions:
\begin{enumerate}
 \item $f^t$ satisfies that $(f^t)'(x)>\sqrt{2}$ for $\forall x\in [-1,1]$.
 \item $f^b_-=\mathrm{id.}$
\end{enumerate}

We depicted the graph of $f^t$ in Figure 9.
Unlike the GLA, the derivative of $f^t$ at $0$ does not diverge to
$\infty$ because of eigenvalues of $\sigma$ and $s$.
However, like the GLA, it is clear that there are infinitely many $f^t$-inverse iteration of 0,
implying $\# (\Lambda \cap l^t)=\infty$.

\begin{figure}[h]
\begin{center}
\includegraphics*[width=5cm]{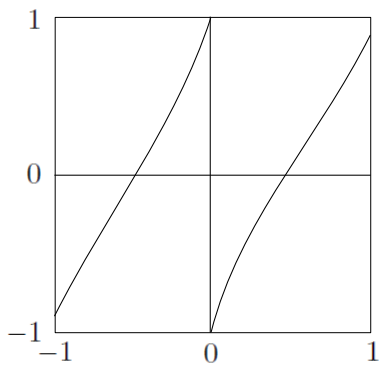}
\end{center}
\caption{1-dimansional map $f^t$}
\end{figure}

By $f^b_-=\mathrm{id.}$, there are infinitely many periodic points in $\Sigma^b_-$.
They accumulate on $W^u_-(\sigma)\subset\Lambda$, hence, $\Lambda$ is not an attracting set.

Now let us see that $\Lambda$ is Lyapunov stable dividing into three parts.
First, since the behaviour of orbits of $\Lambda$ in $\Sigma^t$ is the same as the GLA and
the GLA is attracting, $\Lambda\cap\Sigma^t$ is Lyapunov stable.
Second, in the Cherry-flowbox, for given $U$
we can take a locally positively-invariant neighborhood $V$ (Figure 10).
Third, for given neighborhood $U$ of $W^u_b(s)$,
we can take its neighborhood V such that $X_t(V)\subset U$ ($\forall t\geq 0$) by eigenvalues of $\sigma$ and $s$.
Thus, putting these three parts together, we have checked that $\Lambda$ is Lyapunov stable.\\

\begin{figure}[h]
\begin{center}
\includegraphics*[width=8cm]{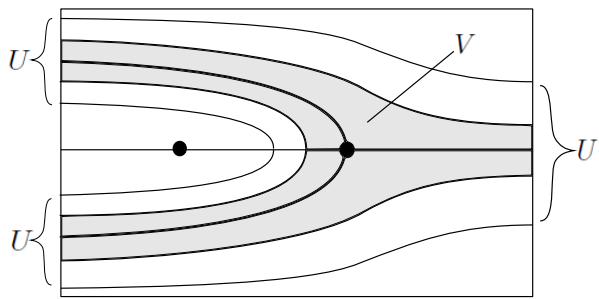}
\caption{Locally positively-invariant neighborhood in the Cherry-flowbox}
\end{center}
\end{figure}

Finally we note that $\Lambda$ does not exhibit the volume-expanding central subbundle, thus $\Lambda$ is not singular-hyperbolic.

%%------------------------------------------------------------------------------------------------------
%%------------------------------------------------------------------------------------------------------

\section*{Acknowledgement}
I am deeply grateful to Prof. S. Hayashi who has been extraordinarily tolerant and gave me insightful comments and suggestions.
Also I received generous supports from Mr. K. Shinohara and Mr. T. Yokoyama.
Finally, I would like to thank Prof. C. Morales for his helpful comments.

%%------------------------------------------------------------------------------------------------------
%%------------------------------------------------------------------------------------------------------

\end{document}